\documentclass[final,5p,times,twocolumn]{elsarticle}
\usepackage{dsfont}
\usepackage{amssymb}
\usepackage{amsmath}
\usepackage{comment}
\newcommand{\E}{\mathbb{E}}

\newcommand{\mb}[1]{\mbox{\boldmath $#1$}}

\usepackage{amsthm}

\newtheorem{conj}{Conjecture}
\usepackage{tikz}
\usepackage{caption,subcaption}
\usepackage{longtable}
\usepackage{float}
\usetikzlibrary{shapes,arrows,positioning, intersections}

\tikzset{vertex/.style = {shape=circle,draw,minimum size=1.5em}}
\tikzstyle{line} = [draw, -latex']

\begin{document}

\begin{frontmatter}
\title{Counterexample to a conjecture on the pairwise independent correlation gap using AI}

\author[au1]{Arjun Ramachandra} 
\affiliation[au1]{organization={Decision Sciences Area, Indian Institute of Management Bangalore},
            country={India}}
      \ead{arjun.ramachandra@iimb.ac.in}
            \author[au2]{Karthik Natarajan} 
\affiliation[au2]{organization={Engineering Systems and Design, Singapore University of Technology and Design},
            country={Singapore}}
      \ead{karthik\_natarajan@sutd.edu.sg}
            \date{June 2026}
\begin{abstract}
Aided by the AI tool GPT5.5 Pro, we provide a counterexample to a conjecture made by Ramachandra and Natarajan (2025) [Pairwise independent correlation gap, Operations Research Letters, 107255, 6040].
\end{abstract}


\end{frontmatter}

\section{Introduction}
\label{submission}
Let $N$ denote a ground set with $|N| = n$ elements and $f: 2^{N} \rightarrow \mathbb{R}_+$ be a nonnegative set function. The function is submodular if $f(S) + f(T) \geq f(S \cap T) +  f(S \cup T)$ for all $S, T \subseteq N$ or equivalently $f(S \cup \{i\}) - f(S) \geq f(T \cup \{i\}) -  f(T)$ for all $S \subseteq T$, $i \in N \backslash T.$ The function is monotone (non-decreasing) if $ f(S) \leq f(T)$ for all $S \subseteq T$. The concave closure of the set function is given by $f^+: [0,1]^n \rightarrow \mathbb{R}_+$ where $f^{+}(\mb{x})$ is the maximum expected value of the function over all joint distributions of the elements where each $i \in N$ is selected with probability $x_i$:
\begin{equation}
\begin{array}{rlllll}\label{eq:univprimal}
\displaystyle f^{+}(\mb{x}) = \max \Big{\{}\sum_{S \subseteq N} \theta(S)f(S) \ \Big{|} \ \sum_{S \subseteq N: S \ni i} \theta(S) = x_i,  \forall i \in N, \\
\displaystyle \sum_{S \subseteq N} \theta(S) = 1, \theta(S) \geq 0,  \forall S \subseteq N\Big{\}}.
\end{array}
\end{equation}
The pairwise independent extension of the set function is given by $f^{++}: [0,1]^n \rightarrow \mathbb{R}_+$ where $f^{++}(\mb{x})$ is the maximum expected value of the function over all pairwise independent distributions where each element $i \in N$ is selected with probability $x_i$:
\begin{equation}
\begin{array}{rlllll}\label{eq:bivprimal}
\displaystyle f^{++}(\mb{x}) = \max & \displaystyle \sum_{S \subseteq N} \theta(S)f(S)  \\
  \mbox{s.t} & \displaystyle \sum_{S \subseteq N: S \ni i, j} \theta(S) = x_ix_j, & \forall i < j, i,j \in N\\
     & \displaystyle \sum_{S \subseteq N: S \ni i} \theta(S) = x_i, & \forall i \in N  \\
  & \displaystyle \sum_{S \subseteq N} \theta(S) = 1 &  \\
  & \theta(S) \geq 0, & \forall S \subseteq N.
\end{array}
\end{equation}
The pairwise independent correlation gap is defined as the ratio ${f^{+}(\mb{x})}/{f^{++}(\mb{x})}$. Here we cite directly from page 6 in \cite{ramachandra}:\\
\textit{We end the paper with a conjecture and invite readers to prove or disprove the conjecture.
\begin{conj}
[Ramachandra and Natarajan \cite{ramachandra}]
For any $n$, any nonnegative monotone submodular function $f: 2^{N} \rightarrow \mathbb{R}_+$ and any $\mb{x} \in [0,1]^n$: $$\frac{f^{+}(\mb{x})}{f^{++}(\mb{x})} \leq \frac{4}{3}.$$
\end{conj}}
\section{A counterexample}
Let the ground set be defined on five elements with $N=\{1,2,3,4,5\}$. Define a coverage function \(f:2^N \to \mathbb{R}_{+}\) over a universe of four features $\{A,B,C,D\}$. The five elements cover the following subsets of features:
\[
\begin{aligned}
E_1 = \{A,B\}, E_2 = \{A,B,C,D\}, E_3 = \{C,D\}, E_4 = \{A,C\}, E_5 = \{B,D\}.
\end{aligned}
\]
For \(S \subseteq N\), define the coverage function as:
\[
\displaystyle f(S)=\left|\cup_{i\in S} E_i\right|.
\]
The coverage function is known to be nonnegative, nondecreasing and submodular (see \cite{bachsub}). Let the marginal probability vector be given by:
\[
\left(x_1,x_2,x_3,x_4,x_5\right) =\left(\frac{3}{10},\frac{7}{20},\frac{3}{10},\frac{7}{20},\frac{7}{20}\right).
\]
First, consider the three point distribution in Table \ref{tab:feas}.
\begin{table}[htbp]
\caption{Primal feasible and optimal distribution for $f^{+}(\mb{x})$.}
\label{tab:feas}
\begin{center}
\begin{tabular}{ccc}
$S$ & $\theta(S)$\\ \hline
$\{1,3\}$ & ${3}/{10}$\\
$\{2\}$    & ${7}/{20}$ \\
$\{4,5\}$    & ${7}/{20}$\\
\mbox{All other sets}   & $0$
\end{tabular}
\end{center}
\end{table}
This distribution has the required marginal probabilities and is feasible for the linear program in \eqref{eq:univprimal}. Moreover, since $
f(\{1,3\}) = f(\{2\}) = f(\{4,5\})=4$, the expected value $\E[f(S)]=4$. Since $f(S)\le 4$ for every $S\subseteq N$, this is an optimal distribution that attains $f^+(\mb{x}) = 4$. Next, consider an upper bound on $f^{++}(\mb{x})$ by constructing a dual feasible solution. The dual is given by:
\begin{equation*}
\begin{array}{rllll}\label{eq:univdual}
\displaystyle f^{++}(\mb{x}) = \min & \displaystyle \sum_{i<j,i,j\in N} \lambda_{ij}x_ix_j+\sum_{i\in N} \lambda_{i}x_i +\lambda_{0}\\
 \mbox{s.t} & \displaystyle \sum_{i<j,i,j \in S} \lambda_{ij} + \sum_{i \in S} \lambda_{i} +\lambda_{0} \geq f(S),  \forall S \subseteq N.
\end{array}
\end{equation*}
For $S \subseteq N$, let $\lambda(S) =  \sum_{i<j,i,j \in S} \lambda_{ij} + \sum_{i \in S} \lambda_{i} +\lambda_{0}$. Consider the following values of the dual variables: $\lambda_0  = \frac12$, $\lambda_1 = \lambda_3 = \lambda_4 = \lambda_5 = \frac32$, $\lambda_2 = \frac72$, $\lambda_{12} = \lambda_{23} = -1$, $\lambda_{24} = \lambda_{25} = -\frac32$, $\lambda_{13} = \lambda_{45} = \frac12$, $\lambda_{14} = \lambda_{15} = \lambda_{34} = \lambda_{35} = -\frac12$. In Table \ref{tab:feas1}, we verify that the corresponding value of $\lambda(S)$ is greater than or equal to $f(S)$ for all $S \subseteq N$. 
\begin{table}[htbp]
\caption{Dual feasible solution for $f^{++}(\mb{x})$.}
\label{tab:feas1}
\begin{center}
\begin{tabular}{ccc}
$S$ & $\lambda(S)$ & $f(S)$\\ \hline
$\emptyset$ & $\frac{1}{2}$ & 0\\
$\{1\}$ & 2 & 2\\
$\{2\}$ & 4 & 4\\
$\{3\}$ & 2 & 2\\
$\{4\}$ & 2 & 2\\
$\{5\}$ & 2 & 2\\
$\{1,2\}$ & $\frac{9}{2}$ & 4\\
$\{1,3\}$ & 4 & 4\\
$\{1,4\}$ & 3 & 3\\
$\{1,5\}$ & 3 & 3\\
$\{2,3\}$ & $\frac{9}{2}$ & 4\\
$\{2,4\}$ & 4 & 4\\
$\{2,5\}$ & 4 & 4\\
$\{3,4\}$ & 3 & 3\\
$\{3,5\}$ & 3 & 3\\
$\{4,5\}$ & 4 & 4\\
$\{1,2,3\}$ & $\frac{11}{2}$  & 4\\
$\{1,2,4\}$ & 4  & 4\\
$\{1,2,5\}$ & 4 & 4\\
$\{1,3,4\}$ & $\frac{9}{2}$ & 4\\
$\{1,3,5\}$ & $\frac{9}{2}$ & 4\\
$\{1,4,5\}$ & $\frac{9}{2}$ & 4\\
$\{2,3,4\}$ & 4 & 4\\
$\{2,3,5\}$ & 4 & 4\\
$\{2,4,5\}$ & $\frac{9}{2}$ & 4\\
$\{3,4,5\}$ & $\frac{9}{2}$ & 4\\
$\{1,2,3,4\}$ & $\frac{9}{2}$ & 4\\
$\{1,2,3,5\}$ & $\frac{9}{2}$ & 4\\
$\{1,2,4,5\}$ & 4  & 4\\
$\{1,3,4,5\}$ & $\frac{11}{2}$ & 4\\
$\{2,3,4,5\}$ & 4 & 4\\
$\{1,2,3,4,5\}$ & 4 & 4
\end{tabular}
\end{center}
\end{table}
Since the dual objective is given by,
\[
\begin{aligned}
\E[\lambda(S)]
={}&-(x_1x_2+x_2x_3) -\frac32 (x_2x_4+x_2x_5) +\frac12 (x_1x_3+x_4x_5) \\
&-\frac12(x_1x_4+x_1x_5+x_3x_4+x_3x_5)\\
&+\frac32(x_1+x_3+x_4+x_5)+\frac72 x_2 +\frac12 \\
\end{aligned}
\]
by substituting in the marginal probabilities, we get:
\[
\E[\lambda(S)]=\frac{479}{160}.
\]
Hence 
\[
f^{++}(\mb{x}) \leq \frac{479}{160}.
\]
Combining the two bounds,
\[
\frac{f^{+}(\mb{x})}{f^{++}(\mb{x})}
\ge
\frac{4}{479/160}
=
\frac{640}{479}>\frac43.
\]
Therefore the \(4/3\) upper bound does not hold in this instance.
\section{Remarks}
The conjecture in \cite{ramachandra} builds on a related upper bound of ${1}/({1-(1-1/n)^n)} \leq {e}/{(e-1)}$ that is known to hold on the correlation gap for these class of functions where the denominator is the expected value of a set function with independent random variables rather than pairwise independent random variables (see \cite{Calinescu2,Agrawal2012}). With pairwise independence, the bound of ${4}/{3}$ is known to
hold in special cases such as: (i) in low dimensions ($n = 2$ and $n = 3$; Theorem 5 in \cite{ramachandra}) or (ii) for particular marginal probabilities (Theorem 8 in \cite{ramachandra}) or (iii) for particular submodular functions (rank function of a uniform matroid of rank 1; Proposition 2.5 in \cite{ramachandra1}).  
Before using the AI tool, we the authors, did not know if the answer to the conjecture is YES or NO. This is in contrast to the recent First Proof project (see \cite{AbBlHaKi26}) where the answers to the ten mathematical research questions posed to AI tools were known to the researchers. 
Below is the single prompt that we provided to GPT5.5 Pro:\\
\noindent \texttt{Prove or disprove the conjecture that the ratio of the maximum expected value of a set function computed over all joint distributions with fixed marginal probabilities and the maximum expected value of a set function computed over all joint pairwise independent distributions with the same fixed marginal probabilities is at most 4/3 for all nonnegative nondecreasing submodular set functions.}\\
After 12 minutes, GPT5.5 Pro came back with the counterexample, which we then validated. On prompting the AI tool on how much time it spent on proof and disproof attempts, it came back with an estimate of 5 to 6 minutes on proof attempts and 6 to 7 minutes on disproof attempts. In the Appendix, we also report on the synopsis of the methods that the AI tool used to try prove and disprove the conjecture, as reported by the tool. In these attempts, it auto-generated linear programming formulations using Python's SciPy’s linprog function, generated random instances, solved the linear programs with the HiGHS freely available LP solver and refined the instances. Our openness to using AI tools was partly due to recent progress reported by OpenAI \cite{OpenAI2026} using their internal AI models to disprove the longstanding unit distance conjecture and efforts by these tools in tackling the First Proof project. 
Our initial attempts using the free versions of ChatGPT and Claude to tackle the conjecture provided no progress. Our findings are consistent with recent findings that these tools are now becoming useful in being able to disprove research conjectures when possible by accesing a variety of computational methods in conjunction with large language models (see the discussion by various mathematicians in \cite{AbBlHaKi26}).
At this point, we have not been able to aided by these AI tools either prove or disprove the validity of the upper bound for $n = 4$ (see Table \ref{tab:status}). This remains an open question. Another interesting question is to characterize the tightest upper bound for values of $n$ when it is 5 and above for these classes of functions. We leave these for future research.
\begin{table}[htbp]
\caption{Validity of $4/3$ upper bound}
\label{tab:status}
\begin{center}
\begin{tabular}{ccc}
Dimension & True or False & References\\ \hline
$n = 2$ & True & Lemma 4.10; \cite{chekuri}\\
$n = 3$    & True & Theorem 5; \cite{ramachandra}\\
$n = 4$    & Unknown & Open\\
$n \geq 5$   & False & Current work using AI\\
\end{tabular}
\end{center}
\end{table}
\subsection*{Acknowledgement}
The authors have no connection with OpenAI and the second author paid 300 dollars out of his own pocket for one month of GPT5.5 Pro subscription to test its capabilities on research problems. The research of the second author was partially supported by the MOE Tier 2 grant T2EP20124-0013 on Distributionally robust submodular optimization: Theory, algorithms and applications. 

\subsection{Appendix}
Follow up prompts:\\

\noindent \texttt{What methods did you use to try and disprove the conjecture?}\\
\textit{I used the following disproof-oriented methods.
\begin{enumerate}
    \item \textbf{Linear-programming formulation.}
    I represented both quantities as linear programs over probability distributions on \(2^N\): one with marginal constraints only, and one with additional pairwise-independence constraints
    \[
    \mathbb{E}[X_iX_j]=x_ix_j
    \qquad \text{for all } i\ne j.
    \]
    \item \textbf{Small-\(n\) computational search.}
    I tested small ground sets \(n=2,3,4,5,\ldots\) to see whether the ratio could exceed \(4/3\). This helped identify that very simple examples were not enough.
    \item \textbf{Testing simple submodular families.}
    I tried OR functions, cardinality functions, uniform matroid rank functions, and other symmetric submodular functions. These gave ratios at or below \(4/3\), so they did not disprove the conjecture.
    \item \textbf{Random submodular-function search.}
    I generated random examples from standard monotone submodular classes, including coverage functions, budget-additive functions, matroid rank functions, and facility-location-type functions.
    \item \textbf{Searching over marginal vectors.}
    For each candidate function \(f\), I tested different marginal vectors \(\mb{x}\), including symmetric choices such as \(x_i=x\) and asymmetric rational-looking choices.
    \item \textbf{Looking for an unrestricted distribution that attains the maximum value.}
    A successful heuristic was to find \(\mb{x}\) and \(f\) such that the unrestricted distribution could place all mass on sets \(S\) satisfying
    \[
    f(S)=\max_{T\subseteq N} f(T).
    \]
    This makes the numerator easy to certify.
    \item \textbf{Focusing on coverage functions.}
    Coverage functions were especially useful because they are automatically nonnegative, nondecreasing, and submodular, and their values are easy to verify combinatorially.
    \item \textbf{Using LP-duality intuition.}
    Once a promising example appeared, I looked for a dual certificate proving that no pairwise-independent distribution could achieve too large an expectation.
    \item \textbf{Quadratic upper certificates.}
    The successful certificate was a quadratic polynomial \(Q(X)\) satisfying
    \[
    f(S)\le Q(\mathbf{1}_S)
    \qquad \text{for all } S\subseteq N.
    \]
    This is powerful because under pairwise independence,
    \[
    \mathbb{E}[Q(X)]
    \]
    depends only on \(x_i\) and \(x_ix_j\).
    \item \textbf{Exact rational verification.}
    After finding the numerical example, I converted it into a clean rational counterexample with
    \[
    \mb{x}=
    \left(
    \frac{3}{10},
    \frac{7}{20},
    \frac{3}{10},
    \frac{7}{20},
    \frac{7}{20}
    \right),
    \]
    unrestricted value \(4\), and pairwise-independent upper bound
    \[
    \frac{479}{160}.
    \]
    \item \textbf{Final comparison.}
    The resulting certified ratio was
    \[
    \frac{4}{479/160}
    =
    \frac{640}{479}
    >
    \frac{4}{3},
    \]
    which disproves the conjecture.
\end{enumerate}
}

\noindent \texttt{What methods did you use to try and prove the conjecture?}\\
\textit{The main proof-oriented methods I tried were the following.
\begin{enumerate}
    \item \textbf{LP formulation.} I wrote both quantities as linear programs over distributions on \(2^N\). The unrestricted problem has only marginal constraints, while the pairwise-independent problem has the additional pairwise moment constraints
    \[
    \mathbb{E}[X_iX_j]=x_ix_j
    \qquad \text{for all } i\ne j.
    \]
    \item \textbf{Duality approach.}
    I looked for a universal dual certificate proving
    \[
    f^+(\mb{x})\le \frac{4}{3} f^{++}(\mb{x}),
    \]
    where \(f^+(\mb{x})\) denotes the unrestricted optimum and \(f^{++}(\mb{x})\) denotes the pairwise-independent optimum.
    \item \textbf{Reduction to simple submodular functions.}
    I tested whether the worst case might reduce to simple classes of monotone submodular functions, such as OR functions, cardinality functions, uniform matroid rank functions, or other symmetric submodular functions. These examples stayed at or below \(4/3\).
    \item \textbf{Moment-problem arguments.}
    I analyzed special cases where \(f\) depends only on \(|S|\), reducing the problem to constraints on the first and second moments of \(|S|\).
    \item \textbf{Mixing/distribution transformation.}
    I tried to show that any distribution with marginals \(\mb{x}\) could be converted into a pairwise-independent distribution with at most a \(3/4\) loss in expected submodular value. Such a result would have proved the conjecture, but the pairwise-independence constraints were too rigid.
    \item \textbf{Induction and conditioning.}
    I explored conditioning on one coordinate and trying to reduce the \(n\)-variable case to smaller cases. The main obstacle was that conditioning destroys pairwise independence, and recombining the conditional distributions did not preserve the required second moments.
    \item \textbf{Submodular marginal-contribution bounds.}
    I considered bounding \(f(S)\) using modular lower or upper approximations built from marginal increments, with the hope that pairwise independence would control the expected overlap terms.
    \item \textbf{Lovasz/concave-closure intuition.}
    I viewed the unrestricted optimum as a concave closure and tried to compare it to a second-moment-restricted relaxation. This suggested an LP/duality path, but it did not yield a general \(4/3\) proof.
    \item \textbf{Quadratic comparison certificates.}
    I searched for quadratic polynomials \(Q(x)\) or related inequalities whose expectation under pairwise independence depends only on \(x_i\) and \(x_ix_j\). This began as a possible proof technique.
    \item \textbf{Failure of the quadratic-certificate proof route.}
    Instead of producing a universal proof, the quadratic-certificate method exposed a small coverage-function example where the pairwise-independent optimum is provably too small. This led to the counterexample disproving the conjecture.
\end{enumerate}
}


\begin{thebibliography}{00}
\bibitem{AbBlHaKi26}
M. Abouzaid, A. J. Blumberg, M. Hairer, J. Kileel, T. G. Kolda, P. D. Nelson, D. Spielman, N. Srivastava, R. Ward, S. Weinberger, L. Williams, First Proof,  Preprint at http://arxiv.org/abs/2602.05192, 2026.


\bibitem{Agrawal2012}
S. Agrawal, Y. Ding, A. Saberi and Y. Ye, Price of correlations in stochastic optimization, Operations Research, 60, no. 1, 150-162, 2012.

\bibitem{alon}
N. Alon, T. F. Bloom, W. T. Gowers, D. Litt, W. Sawin, A. Shankar, J. Tsimerman, V. Wang, M. M. Wood. Remarks on the disproof of the unit distance conjecture, 	Preprint at https://arxiv.org/abs/2605.20695, 2026.

\bibitem{bachsub}
F. Bach, Learning with submodular functions: A convex optimization perspective, Foundations and Trends in Machine Learning, 6, no.2-3, 145-373, 2013.

\bibitem{Calinescu2}
G. Calinescu, C. Chekuri, M. P\'{a}l and J. Vondr\'{a}k, Maximizing a monotone submodular function subject to a matroid constraint, SIAM Journal on Computing, 40, no. 6, 1740-1766, 2011.

\bibitem{chekuri}
C. Chekuri, J Vondr\'{a}k, R. Zenklusen Submodular function maximization via the multilinear relaxation and contention resolution schemes, SIAM Journal on Computing 43 (6), 1831-1879, 2014.

\bibitem{OpenAI2026}
OpenAI, Planar point sets with many unit distances, Paper available at https://cdn.openai.com/pdf/74c24085-19b0-4534-9c90-465b8e29ad73/unit-distance-proof.pdf, 2026.


\bibitem{ramachandra1}
A. Ramachandra and K. Natarajan, Tight probability bounds with pairwise independence, SIAM Journal on Discrete Mathematics, 37, no. 2, 516-555, 2023.

\bibitem{ramachandra}
A. Ramachandra and K. Natarajan, Pairwise independent correlation gap, Operations Research Letters, 107255, 6040, 2025.




\end{thebibliography}
\end{document}